\documentclass[12pt]{article}
\usepackage{amsfonts,amsmath,amssymb,amscd}

\begin{document}

\parskip 1pc

\renewcommand{\thefootnote}{\arabic{footnote}}
\newcommand{\bea}{\begin{eqnarray}}
\newcommand{\ena}{\end{eqnarray}}
\newcommand{\beas}{\begin{eqnarray*}}
\newcommand{\enas}{\end{eqnarray*}}
\newtheorem{theorem}{Theorem}[section]
\newtheorem{corollary}{Corollary}[section]
\newtheorem{conjecture}{Conjecture}[section]
\newtheorem{proposition}{Proposition}[section]
\newtheorem{lemma}{Lemma}[section]
\newtheorem{definition}{Definition}[section]
\newtheorem{example}{Example}[section]
\newtheorem{remark}{Remark}[section]

\def\proof{{\bf Proof}}
\def\qed{\hfill \mbox{\rule{0.5em}{0.5em}}}
\newcommand\ape[1]{}
\newcommand{\ignore}[1]{}
\def\bbbone{{\mathchoice {\rm 1\mskip-4mu l}
{\rm 1\mskip-4mu l}{\rm 1\mskip-4.5mu l} {\rm 1\mskip-5mu l}}}
\def\proof{\noindent{\bf Proof.\ \ }}

\def\todist{\stackrel{\mbox{distr}}{\longrightarrow}}
\def\indist{\stackrel{\mbox{\tiny d}}{=}}

\def\bbbe{{\rm I\!E}}
\def\bbbp{{\rm I\!P}}
\def\bbbr{{\rm I\!R}}
\def\p{\bbbp}   
\def \e {\bbbe}   
\def\BR{\bbbr}
\def\BZ{\bbbz}
\def\levy{L\'evy }

\title
{Size bias, sampling, the waiting time paradox, and infinite divisibility: when is the increment independent?}

\author{Richard Arratia, Larry Goldstein}


\maketitle

\abstract
With $X^*$ denoting a random variable with the $X$-size bias distribution,
what are  all distributions for $X$ such that
it is possible to have $X^*=X+Y$, $Y\geq 0$,
with $X$ and $Y$ {\em  independent}?  We give the answer, due to Steutel \cite{steutel},
and also discuss the relations of size biasing to
the waiting time paradox, renewal theory, sampling, tightness and uniform integrability,
compound Poisson
distributions, infinite divisibility, and the lognormal distributions.
 \endabstract


\section{The Waiting Time Paradox}
\label{sec:waiting}

Here is the ``waiting time paradox,'' paraphrased from
Feller \cite{feller2}, volume II, section I.4:
Buses arrive in accordance with a Poisson process, so that
the interarrival times are
given by  independent random variables, having the
exponential distribution $\p(X >s)=e^{-s}$ for $s>0$, with mean $\e X = 1$.  I now arrive at an arbitrary time $t$.  What is the expectation
$\e W_t$ of my waiting time $W_t$ for the next bus?  Two contradictory answers
stand to reason: (a) The lack of memory of the exponential distribution, i.e.~the property $\p(X>r+s|X>s)=\p(X>r)$, implies that $\e W_t$ should not be sensitive to the
choice $t$, so that $\e W_t = \e W_0 =1$.  \ (b) The time of my arrival is
``chosen at random'' in the interval between two consecutive buses, and for reasons
of symmetry $\e W_t = 1/2$.

The resolution of this paradox requires an understanding of
size biasing.   We will first present some simpler examples of size biasing, before
returning to
the waiting time paradox and its resolution.

Size biasing occurs in many unexpected contexts, such as statistical estimation, renewal theory, infinite divisibility of distributions, and number theory.
The key relation is that  to size bias a sum with
independent summands, one needs only size bias a single summand, chosen at
random.

\section{Size Biasing in Sampling}
\label{sect sampling}

We asked  students who ate lunch in the cafeteria ``How many people, including yourself,  sat at your table?''  Twenty percent said they ate alone, thirty percent said they ate with one other person, thirty percent said they ate at a table of three, and the remaining twenty percent said they ate at a table of four.  From this
information, would it be correct to conclude  that twenty percent of the   tables had only one person,
thirty percent had two people,
thirty percent had three people, and
twenty percent had four people?

Certainly not!  The easiest way to think about this situation is to imagine
100 students went to lunch, and we interviewed them all.  Thus, twenty students ate
alone, using 20 tables, thirty students ate in pairs, using 15 tables, thirty
students ate in trios, using 10 tables, and twenty students ate in groups of four,
using 5 tables.  So there were $20+15+10+5=50$ occupied tables, of which
forty percent had only one person, thirty percent had two people, twenty percent had three people, and ten percent had four people.

A probabilistic view of this example begins by considering the experiment where an occupied table is selected at random and the number of people, $X$, at that table is recorded. From the analysis so far, we see that since 20 of the 50 occupied tables had only a single individual, $\p(X=1)=.4$, and so forth. A different experiment, one related to but not to be confused with the first, would be to select a person at random, and record the total number $X^*$ at the table where this individual had lunch. 
Our story \emph{began} with the information $\p(X^*=1)=.2$, $\p(X^*=2)=.3$, and so forth, and the 
distributions of the random variables $X$ and $X^*$ are given side by side in the following 
table:
$$
\begin{array}{|c|c|c|} \hline
k & \p(X=k) &  \p(X^*=k) \\ \hline
1 & .4 & .2 \\
2 & .3 & .3 \\
3 & .2 & .3 \\
4 & .1 & .2 \\ \hline
 & 1.0 & 1.0 \\ \hline
 \end{array}
$$
The distributions of the random variables $X$ and $X^*$ are related;
for  $X$ each  table  has  the  same  chance  to  be  selected, but  for
$X^*$  the  chance  to  select  a  table  is  proportional  to  the  number  of
people  who  sat  there.  Thus  $\p(X^*=k)$ is  proportional  to  $k \times
\p(X=k)$;  expressing  the  proportionality  with  a constant $c$ we have
\mbox{$\p(X^*=k)$}  $= c \times \p(X=k)$.  Since $1=\sum_k \p(X^*=k) = c \sum_k
k \p(X=k)=c \e X$, we have $c=1/\e X$ and
\begin{equation}\label{bias-k}
\p(X^*=k) = \frac{k \p(X=k)}{\e X};   \ \ k=0,1,2,\ldots.
\end{equation}
Since the distribution of $X^*$ is weighted by the value, or size, of $X$, we say that $X^*$ has the $X$ \emph{size biased distribution}.

In many statistical sampling situations, like the one above, care must be taken so that one does not inadvertently sample from the size biased distribution in place of the one intended. For instance, suppose we wanted to have information on how many voice telephone lines are connected at residential addresses. Calling residential telephone numbers by random digit dialing and asking how many telephone lines are connected at the locations which respond is an instance where one would be observing the size biased distribution instead of the one desired. It's three times more likely for a residence with three lines to be called than a residence with only one. And the size bias distribution never has any mass at zero, so no one answers the phone and tells a surveyor that there are no lines at the address just reached! But the same bias exists more subtly in other types of sampling more akin to the one above: what if we were to ask people at random how many brothers and sisters they have, or how many fellow passengers just arrived with them on their flight from New York?

\section{Size Bias in General}

The  examples in Section \ref{sect sampling} involved nonnegative integer
valued random variables.  In general, a random variable $X$ can be size biased if and only if it is nonnegative, with finite and positive mean, i.e. $1=\p(X \ge 0)$ and $0< \e X < \infty$.   We will henceforth
assume that $X$ is nonnegative, with  $a := \e X \in (0,\infty)$. For such $X$, we say $X^*$ has the $X$ size
biased distribution if and only for all
bounded continuous functions $g$,
\begin{equation}\label{sizebias}
\e  g(X^*) =  \frac{1}{a} \ \e (X g(X) ).
\end{equation}
It is easy to see that, as a condition on distributions, (\ref{sizebias}) is equivalent to
\beas
dF_{X^*}(x) = \frac{  x \, dF(x)}{a}.
\enas
In particular, when $X$ is discrete with probability mass function $f$, or when $X$ is continuous
with density $f$, the formula
\bea \label{sizebias-f}
f(x)=\frac{xf(x)}{a},
\ena
applies; \eqref{bias-k} is a special case of the former.

If \eqref{sizebias} holds for all bounded continuous $g$, then by
monotone convergence it also holds for
any function $g$ such that $\e |X g(X)| < \infty$.  In particular,
taking
  $g(x)=x^n$, we have
\begin{equation}\label{moment shift}
\e (X^*)^n=\e X^{n+1}/\e X
\end{equation}
whenever $\e |X^{n+1}|<\infty$.  Apart from the extra scaling by $1/\e X$, \eqref{moment shift} says that the sequence of moments
of $X^*$ is the sequence of moments of $X$, but shifted by one.
One  way to recognize size biasing is through
the ``shift of the moment sequence;'' we give an example in Section
 \ref{sect lognormal}.

In this paper, we ask and solve the following problem:  what are all possible distributions for $X \geq 0$ with $0< \e X<\infty$, such
that there exists a coupling in which
\begin{equation}\label{indep}
X^*=X+Y, \ \ Y \geq 0, \ \ \ \mbox{ and } X,Y \mbox{ are independent}.
  \end{equation}
Resolving this question on independence leads us to the
infinite divisible and compound Poisson distributions.  These concepts by themselves
can be quite technical, but in our size biasing context they are relatively easy.
We also present some background information on size biasing,
in particular how it arises in applications including statistics.
The answer to (\ref{indep}) comes from Steutel 1973 \cite{steutel}; see section 10 for
more of the history.

A beautiful treatment of
size biasing for branching processes is \cite{LPP} by Lyons, Pemantle, and Peres.
Size biasing has a connection with Stein's method for obtaining error bounds when
approximating the distributions of sums by the Normal, (\cite{BR89} Baldi, P. Rinott, Y. 1989, \cite{BRS89} Baldi, P. Rinott, Y.  and Stein C., 1989, and
\cite{GR96} Goldstein and Rinott, 1996), and the Poisson
(\cite{BHJ92}, Barbour, Holst, and Janson, 1992).

To more fully explain the term ``increment'' in the title,
letting $g(x) = \bbbone(x>t)$ in (\ref{sizebias}) for some fixed $t$, we find that
\[
\p(X^*>t)= \frac{1}{a}\ \e (X\bbbone(X>t)) \ \ge \  \frac{1}{a}\ \e X \  \e  \bbbone(X>t) =\p(X>t).
\]
The inequality above is the special case $f(x)=x$, $g(x) =\bbbone(x>t) $
of Chebyschev's correlation inequality: $\e (f(X) g(X)) \geq
\e f(X) \ \e g(X)$ for any random variable and any
two increasing functions $f,g$.
The condition $\p(X^*>t) \geq \p(X>t)$ for all $t$ is described as ``$X^*$ lies
above $X$ in distribution,'' and implies  that there exist couplings of
$X^*$ and $X$ in which always $X^* \geq X$.  Writing $Y$ for  the
difference, we have
 \begin{equation}
\label{xplusy}
X^*=X+Y, \ \ Y \geq 0.
  \end{equation}

The simplest coupling satisfying (\ref{xplusy})
is based on the ``quantile transformation,''
constructing each of $X$ and $X^*$ from the same uniform
random variable $U$ on (0,1). Explicitly, with cumulative distribution
function $F$ defined by $F(t):= \p(X \leq t)$, and its ``inverse'' defined by
$F^{-1}(u) := \sup \{ t\!: \ F(t) \leq u \} $, the coupling given by
$X=F^{-1}(U), X^*=(F^*)^{-1}(U)$ satisfies (\ref{xplusy}).

In general (\ref{xplusy})
determines neither the joint distribution of $X$ and $Y$, nor  the
marginal distribution of $Y$, nor whether or not $X$ and $Y$ are independent.
It is a further restriction on the distribution of $X$
 to require that (\ref{xplusy}) be achievable with
$X,Y$ {  \em independent}.

When $Z \sim Po(\lambda)$, i.e. $Z$ is Poisson with
$\p(Z=k)=e^{-\lambda} \lambda^k/k!, \ $ $k=0,1,2,\ldots$,
we have $Z^* \indist Z+1$,
where the notation $\indist$ denotes equality in distribution.
The reader can check
\begin{equation} \label{monkey1}
Z^* \indist Z+1
\end{equation}
 directly using
(\ref{sizebias}); a conceptual derivation is given in Example 1)
in Section \ref{exnoindep}.
Scaling  by a factor $y>0$ in general means to replace $X$ by $yX$, and it follows
easily from (\ref{sizebias}) that
\begin{equation}\label{scaling}
(yX)^*=y(X^*).
\end{equation}
For our case, multiplying (\ref{monkey1}) by $y > 0$ yields
the implication, for Poisson $Z$,
\begin{equation}
\label{yXpois}
\mbox{if} \quad X = yZ,   \quad \mbox{then} \quad  X^*=X+y.
\end{equation}

Hence, for each $\lambda>0$ and $y>0$, (\ref{yXpois}) gives an example where (\ref{indep}) is satisfied
with $Y$ a constant
random variable, which is independent of {\em every}
random variable.  In a very concrete sense, all solutions of (\ref{indep})
can be built up from these examples, but to accomplish that we must
first review how to size bias sums of independent random variables.

\section{How to size bias a sum of independent random variables}
Consider a sum $X=X_1+\cdots+X_n$,
with independent non-negative summands $X_i$,
and suppose that $ \e X_i=a_i$, $ \e X=a$.
Write $S_i=X-X_i$, so that $S_i$ and $X_i$ are independent, and
also take  $S_i$  and $X_i^*$  to be independent;  this is used to obtain
the final inequality in (\ref{mixture}) below.

We have for all bounded functions $g$,

\bea
  \e g(X^*) &=&  \e (Xg(X))/a\nonumber \\
         &=& \sum_{i=1}^n  (a_i/a) \e (X_i g(S_i+X_i))/a_i \nonumber \\
         &=& \sum_{i=1}^n (a_i/a) \e g(S_i+X_i^*). \label{mixture}
  \ena
The result in (\ref{mixture}) says precisely
that $X^*$ can be represented by the mixture of the distributions
$S_i+X_i^*$ with mixture probabilities $a_i/a$. In  words, in order to
size bias the sum $X$ with independent summands,
we first pick an independent index $I$ with probability proportional to
its expectation, that is, with distribution $\p(I=i)=a_i/a$, and then size bias only the summand $X_I$.
Or, with $X_1,\ldots,X_n,X_1^*,\ldots,X_n^*$ and $I$ all independent
\begin{equation}
\label{noniid}
(X_1+X_2+\cdots+X_n)^*=X_1+ \cdots + X_{I-1}+X_I^*+X_{I+1}+\cdots+X_n.
\end{equation}

For the special case where the summands $X_i$ are not only independent but also {\em identically distributed}, or i.i.d., this recipe simplifies.
In this case it does not matter which summand is biased, as all the distributions in the
mixture are the same; hence for any $i=1,\ldots,n$, $X^* \indist X_1+\cdots+
X_{i-1}+X_i^*+X_{i+1}+\cdots+X_n$.  In particular we may use $i=1$ so that
\begin{equation}\label{iid}
(X_1+X_2+\cdots+X_n)^*=X_1^*+ X_{2}+X_3+\cdots+X_n.
\end{equation}

\section{Waiting for a bus: the renewal theory  connection}\label{sec:renewal}
Renewal theory provides a
conceptual explanation of the identity
(\ref{iid}) and at the same time gives an explanation of the waiting time paradox.
Let the interarrival times of our buses in Section \ref{sec:waiting} be denoted $X_i$,
so that buses arrive at
times $X_1, X_1+X_2, X_1+X_2+X_3,\ldots$,
and assume only that the $X_i$ are
i.i.d., strictly positive random variables with
finite mean; the  paradox presented earlier was the special case with $X_i$ exponentially distributed.
Implicit in the story of my arrival time $T$ as ``arbitrary'' is that my  precise arrival time does not matter, and that there should be
no relation
between my arrival time and the schedule of buses.  One
way to model this assumption is to choose $T$ uniformly from 0 to $l$, independent of
$X_1,X_2,\ldots$, and then take
the limit as $l \rightarrow \infty$; informally, just imagine
some very large $l$. Such a $T$ corresponds
to throwing a dart at random from a great distance toward the real line, which has
been subdivided into intervals of lengths $X_i$.  Naturally the dart is twice as likely to land
in a given interval of length two than one of length one, and generally $x$ times as likely to land in a given
interval of length $x$ as one of length one. In other words, if the interarrival times $X_i$ have a distribution $dF(x)$, the distribution of the length of the interval where the dart lands is
proportional to $x \ dF(x)$.  The constant of proportionality must be $1/a$,
in order to make a legitimate distribution,
so the distribution of the interval where the dart lands is the distribution of $X^*$.

The conceptual explanation of identity (\ref{iid}) is the following.
Suppose that every $n^{th}$ bus is bright blue, so that the waiting time between
bright blue buses is the sum over a block of $n$ successive arrival times.
Again, the random time $T$ finds itself in an interval whose length is
distributed as the size biased
distribution of the interarrival times; the length of the neighboring intervals are not affected.
But by considering the variables as appearing in blocks of $n$, the random time $T$ must also find
itself in a block distributed as $(X_1+\cdots+X_n)^*$. Since only
the interval containing one of the interarrival times has been size biased, this sum must
be equal in distribution to $X_1+\cdots+X_{i-1}+X_i^*+X_{i+1}+\cdots+X_n$.

A more precise explanation of our waiting time paradox is based on the concept of
stationarity ---
randomizing the schedule of buses so that I can
 arrive at an arbitrary time $t$,
and specifying a particular $t$ does not
influence how long I must wait for the next bus.
The simple process  with arrivals at times $X_1, X_1+X_2, X_1+X_2+X_3,\ldots$  is in general not
stationary;  and the distribution of the time $W_t$  that we wait from time $t$ until the arrival of
the next bus  varies with $t$.  We can, however, cook up a stationary process from this simple process by a modification suggested by size biasing.
For motivation, recall the case where I  arrive at $T$ chosen
uniformly from $(0,l)$.
In the limit as $ l \rightarrow \infty$ the interval containing $T$ has length distributed as
$X_i^*$, and my arrival within this
interval is `completely random.'  That is,
I wait  $U X_i^*$ for the next bus,
and I missed the previous bus by $(1-U)X_i^*$, where $U$ is uniform on (0,1) and independent of $X_i^*$.
Thus it is plausible that one
can form a stationary renewal process by the following recipe.
Extend $X_1,X_2,\ldots$ to an independent, identically distributed sequence
$\ldots,X_{-2},X_{-1},X_0,X_1,X_2,\ldots$ \ . Let $X_0^*$ be
the size biased version of $X_0$  and let $U$ be chosen uniformly in (0,1), with
all variables independent.
The origin is to occupy an
interval of length $X_0^*$, and the location of the origin is to be uniformly distributed over this interval; hence buses arrive at time $UX_0^*$ and $-(1-U)X_0^*$.
Using $X_1,X_2,\ldots$ and $X_{-1},X_{-2},\ldots$ as interarrival times on the
positive and negative side, we obtain a process by
setting bus arrivals
at the positive times $UX_0^*, UX_0^*+X_1,UX_0^*+X_1+X_2,\cdots$, and at the negative times
$-(1-U)X_0^*,-((1-U)X_0^*+X_{-1}),-((1-U)X_0^*+X_{-1}+X_{-2}),\ldots$\ \ ,
and it can be proved that this process is stationary.

The interval which covers the origin has expected length
$\e X_0^* = \e X_0^2 / \e X_0$ (by (\ref{moment shift}) with $n=1$,)
and the ratio of this to $\e X_0$ is $\e X_0^*/\e X_0 = \e X_0^2 / (\e X_0)^2$.
By Cauchy-Schwarz, this ratio is at least 1;  and every value in $[1,\infty]$ is feasible.
Note that my waiting time is $\e W_T = \e W_0 = \e (U X_0^*) = (1/2) \e X_0^*$,
so the ratio of my waiting time to the average time between buses can be any value
between 1/2 and infinity, depending on the distribution of the interarrival times.

The exponential case is very special, where strange and wonderful
``coincidences''
effectively hide all the structure involved in size biasing and stationarity.
The distribution of $X_0^*$, obtained by size biasing the unit exponential,
has density $xe^{-x}$ for $x>0$, using (\ref{sizebias-f}) with $a=1$.
This distribution is known as Gamma(1,2).
In particular, $\e X_i^* = \int_0^\infty x (xe^{-x}) \ dx =2$,
and splitting this in half for ``symmetry'' as in Feller's answer (b) gives
1 as the expected time I must wait for the next bus.
Furthermore, the independent uniform $U$ splits that Gamma(1,2) variable $X_0^*$ into $UX_0^*$ and
$(1-U)X_0^*$, and these turn out to be independent, and each having the original exponential distribution. Thus the general recipe for cooking up a stationary
process, involving $X_0^*$ and $U$ in general, simplifies beyond recognition:
the original simple schedule with arrivals at times $X_1, X_1+X_2,
X_1+X_2+X_3,\ldots$ forms half of a stationary process,
which is completed by its other half, arrivals at $-X_1', -(X_1'+X_2'), \ldots, $
with $X_1,X_2,\ldots,X_1',X_2',\ldots$ all independent and
exponentially distributed.

\section{Size bias in statistics}

But size biasing is not always undesired. In fact, it can be used to construct
{\em unbiased} estimators of quantities that are at first glance difficult to
estimate without bias. Suppose we have a population of $n$ individuals, and
associated to each individual $i$ is the pair of real numbers $x_i \ge 0$ and $y_i$,
with $\sum x_i >0$. Perhaps $x_i$
is how much the $i^{th}$
customer was billed by their utility company last month, and $y_i$, say a smaller
value than $x_i$, the amount they were supposed to have been billed.
Suppose we would like to know just how severe the overbilling error is; we would like
to know the
'adjustment factor', which is the ratio $\sum_i y_i / \sum_i x_i$. Collecting
the paired values for everyone is laborious and expensive, so we would like to be able to
use a sample of $m<n$ pairs to make an estimate. It is not too hard to verify that
if we choose a set $R$ by selecting $m$ pairs uniformly from the $n$, then the estimate
$\sum_{j \in R}y_j/\sum_{j \in R}x_j$ will be biased; that is, the estimate, on average,
will not equal the ratio we are trying to estimate.

 Here's how size biasing can be used to construct an unbiased estimate of the ratio $\sum_i y_i / \sum_i x_i$, using $m<n$ pairs. Create a random set $\tilde{R}$ of size $m$ by first selecting a pair with probability proportional to $x_i$, and then $m-1$ pairs uniformly from the remaining pairs. Though we are out of the independent framework, the principle of (\ref{iid}) is still at work; size biasing one has size biased the sum. (This is so because we have size biased the
one, and then chosen the others from an appropriate conditional distribution.) That is,
one can now show that by biasing to
include the single element in proportion to its $x$ value, we have achieved a distribution whereby the
probability of choosing the set $r$ is proportional to $\sum_{j \in r}x_j$. From this observation it
is not hard to see why $\e(\sum_{j \in \tilde{R}}y_j/\sum_{j \in \tilde{R}}x_j)=\sum_i y_i / \sum_i x_i$. This method
is known as Midzuno's procedure for unbiased ratio estimation, and is noted in
Cochran  \cite{C77}.

\section{Size biasing, tightness, and uniform integrability}

Recall that a collection of random variables $\{Y_\alpha: \alpha \in I \}$ is \emph{tight} iff for all $\varepsilon>0$ there exists $L<\infty$ such that
$$
\p(Y_\alpha \not \in [-L,L])  < \varepsilon \quad \mbox{for all $\alpha \in I$.}
$$
This definition looks quite similar to the definition of uniform integrability, where we say
$\{X_\alpha: \alpha \in I \}$ is \emph{uniformly integrable}, or UI, iff for all $\delta>0$ there exists $L<\infty$ such that $$
\e (|X_\alpha| ; X_\alpha \notin [-L,L])  < \delta \quad \mbox{for all $\alpha \in I$.}
$$
Intuitively, tightness for a family is that uniformly over the family, the probability mass due to large values
is arbitrarily small.  Similarly, uniform integrability is the condition that, uniformly over the family, the contribution to the expectation due to large values is arbitrarily small.

Tightness of the family of random variables $\{Y_\alpha: \alpha \in I \}$ implies that every sequence of variables $Y_n,n=1,2,\ldots$ from the family has a subsequence that converges in distribution. The concept of tightness is very useful not just for random variables, that is,  real-valued random  objects,
but also for random elements of other spaces; in more general spaces, the closed intervals $[-L,L]$ are replaced by
 \emph{compact sets}. If $\{X_\alpha: \alpha \in I\}$ is uniformly integrable, $\e X_n \rightarrow \e X$ for any sequence of variables $X_n, n=1,2,\ldots$ from the family that converges in distribution.

To discuss the connection between size biasing and uniform integrability, it is useful to restate the basic definitions in terms of nonnegative random variables. It is clear from the definition of tightness above that a family of \emph{nonnegative} random variables $\{Y_\alpha: \alpha \in I \}$ is tight iff for all $\varepsilon>0$ there exists $L<\infty$ such that
\begin{equation}\label{def tight}
\p(Y_\alpha > L)  < \varepsilon  \quad \mbox{for all $\alpha \in I$,}
\end{equation}
and from the definition of UI, that
a family of \emph{nonnegative} random variables  $\{X_\alpha: \alpha \in I \}$ is uniformly integrable iff for all $\delta>0$ there exists $L<\infty$ such that
\begin{equation}\label{def UI}
\e(X_\alpha ; X_\alpha > L)  < \delta \quad \mbox{for all $\alpha \in I$.}
\end{equation}
For general random variables, the family  $\{G_\alpha: \alpha \in I \}$ is tight [respectively UI]
iff $\{|G_\alpha|: \alpha \in I \}$ is tight [respectively UI].
We specialize in the remainder of this section to random variables that are non-negative with finite, strictly positive  mean.

Since \emph{size bias} relates contribution to the expectation to probability mass, there should be a connection between tightness, size bias, and UI.  However, care should be taken to distinguish between the (additive) contribution to expectation, and the \emph{relative} contribution to expectation. The following example makes this distinction clear. Let
$$
\p(X_n= n)=1/n^2, \p(X_n = 0) = 1 - 1/n^2, \quad n=1,2,\ldots.
$$
Here, $\e X_n = 1/n$, the family $\{X_n\}$ is uniformly integrable, but $1=\p(X_n^*=n)$, so the family
$\{X_n^*\}$ is not tight. The trouble is that the additive contribution to the expectation from large values of $X_n$ is small, but the \emph{relative} contribution is large --- one hundred percent!
The following two theorems, which exclude this phenomenon, show that tightness and uniform integrability are very closely related.

\begin{theorem}
\label{thm:tui1}
Assume that for $\alpha \in I$, where $I$ is an arbitrary index set, the random variables $X_\alpha$ satisfy
$X_\alpha \ge 0$ and $c \le \e X_\alpha < \infty$, for some $c>0$.  For each $\alpha$ let $Y_\alpha = X_\alpha^*$.
 Then
$$
\{ X_\alpha: \alpha \in I \} \mbox{ is UI \ iff } \  \{ Y_\alpha: \alpha \in I \} \mbox{ is tight}.
$$
\end{theorem}

\proof First, with $Y_\alpha = X_\alpha^*$, we have  $\p(Y_\alpha > L ) = \e ( 1(Y_\alpha > L) ) =
\e (X_\alpha 1(X_\alpha >L)) / \e X_\alpha$,  so for any $L$ and $\alpha \in I$,
$$
\e(X_\alpha ;X_\alpha >L) = \e X_\alpha \p(Y_\alpha > L ).
$$

Assume that $\{ X_\alpha: \alpha \in I \}$ is UI, and let $\varepsilon >0 $ be given to test tightness in \eqref{def tight}. Let $L$ be such that \eqref{def UI} is satisfied with $\delta = \varepsilon c$.  Now, using
$\e X_\alpha \ge c$, for every $\alpha \in I$,
$$
 \p(Y_\alpha > L ) = \e(X_\alpha ;X_\alpha >L) / \e X_\alpha \le \e(X_\alpha ;X_\alpha >L) /c < \delta /c =\varepsilon,
$$
establishing \eqref{def tight}.

Second, assume that $\{ X_\alpha: \alpha \in I \}$ if tight, and take  $L_0$ to satisfy \eqref{def tight} with $\varepsilon := 1/2 $, so that $\p(Y_\alpha > L_0 ) <1/2$ for all $\alpha \in I$.  Hence, for all $\alpha \in I$,
$$
\e(X_\alpha ;X_\alpha >L_0) = \e X_\alpha  \p(Y_\alpha > L_0 ) <  \e X_\alpha / 2,
$$
and therefore,
\beas
L_0 \ge  \e(X_\alpha ;X_\alpha  \le L_0) &=& \e X_\alpha - \e(X_\alpha ;X_\alpha >L_0)\\
&>& \e X_\alpha - \e X_\alpha / 2 = \e X_\alpha / 2,
\enas
and hence $\e X_\alpha < 2 L_0$.
Now given $\delta >0$ let $L$ satisfy \eqref{def tight} for $\varepsilon= \delta/(2 L_0)$.
Hence $\forall \alpha \in I$,
$$
  \e(X_\alpha ;X_\alpha >L) = \e X_\alpha \ \p(Y_\alpha > L ) < 2 L_0 \  \p(Y_\alpha > L ) < 2 L_0 \ \varepsilon = \delta,
$$
establishing \eqref{def UI}.
\\

\begin{theorem}
\label{thm:tui2}
Assume the for $\alpha \in I$, where $I$ is an arbitrary index set, that random variables $X_\alpha$ satisfy
$X_\alpha \ge 0$ and $\e X_\alpha < \infty$.  Pick any $c \in (0,\infty)$, and for each $\alpha$ let $Y_\alpha = (c+X_\alpha)^*$.
Then
$$
\{ X_\alpha: \alpha \in I \} \mbox{ is UI \ iff } \  \{ Y_\alpha: \alpha \in I \} \mbox{ is tight}.
$$
\end{theorem}

\proof  By Theorem \ref{thm:tui1}, the family $\{c+X_\alpha\}$ is UI iff the family $\{(c+X_\alpha)^*\}$ is tight.
As it is easy to verify that the family $\{X_\alpha\}$ is tight [respectively UI] iff the family $\{c+X_\alpha\}$ is tight [respectively UI], Theorem \ref{thm:tui2} follows directly from Theorem \ref{thm:tui1}.

\section{Size biasing and infinite divisibility: the heuristic}\label{sect inf div}
Because of the recipe (\ref{iid}), it is natural that our question
in (\ref{indep}) is related to the
concept of infinite divisibility. We say that a random variable $X$ is infinitely
divisible if for all $n$, $X$ can be decomposed in distribution
as the sum of $n$ iid variables. That is, that for all $n$ there exists a distribution
$dF_n$ such that if $X_1^{(n)},\ldots,X_n^{(n)}$ are iid with this distribution, then
\begin{equation}\label{inf div}
X \indist X_1^{(n)}+\cdots + X_n^{(n)}.
\end{equation}
Because this is an iid sum, by (\ref{iid}), we have
$$
X^*=(X-X_1^{(n)}) \ + \ (X_1^{(n)})^*,
$$
with $X-X_1^{(n)}$ and $(X_1^{(n)})^*$ independent.
For large $n$, $X-X_1^{(n)}$ will be close to $X$, and so we have represented the size
bias distribution of $X$ as approximately equal, in distribution, to $X$ plus an
independent increment. Hence it is natural to suspect that the class of non negative
infinitely divisible random variables can be size biased by adding an independent
increment.

It is not difficult to make the above argument rigorous, for
infinitely divisible $X \ge 0$ with
$\e X < \infty$.
First, to show that $X-X_1^{(n)}$
converges in distribution to $X$, it suffices
to show that $X_1^{(n)}$ converges to zero
in probability.
Note that $X \ge 0$ implies $X_1^{(n)} \ge 0$, since \eqref{inf div}
gives
$0 = \p(X<0) \ge (\p(X_1^{(n)}<0)^n$.
Then, given $\epsilon >0$, $\infty > \e X \ge n \p(X_1^{(n)} >
\epsilon) \epsilon$  implies that
$\p(X_1^{(n)}>\epsilon) \to 0$; hence  $X_1^{(n)} \rightarrow 0$
in probability as $n \rightarrow \infty$.

We have that
\begin{equation}\label{inf div sum}
X^* = (X-X_1^{(n)}) + (X_1^{(n)})^*,
\end{equation}
with $X-X_1^{(n)}$ and $(X_1^{(n)})^*$ independent,
and $X-X_1^{(n)}$ converging to $X$ in distribution.
Now, the family of random variables $(X_1^{(n)})^*$ is ``tight'', because given
$\epsilon >0$, there is a $K$ such that $\p(X^* >K)<\epsilon$, and by
(\ref{inf div sum}), for
all $n$, $\p((X_1^{(n)})^*>K) \leq \p(X^*>K) <\epsilon$.  Thus, by Helly's theorem,
there exists a
subsequence $n_k$ of the $n$'s along which $(X_1^{(n)})^*$ converges in distribution,
say $(X_1^{(n_k)})^* \todist Y$.
Taking  $n \rightarrow \infty$ along this subsequence, the pair  $(X-X_1^n,(X_1^n)^*)$
converges jointly to the pair $(X,Y)$ with $X$ and $Y$ independent.  From
$X^* = (X-X_1^{(n_k)}) + (X_1^{(n_k)})^* \todist X+Y$ as $k \rightarrow \infty$
we conclude that $X^* \indist X+Y$, with $Y \geq 0$, and  $X,Y$ independent.
This concludes a proof that if $X\geq 0$ with $0<\e X < \infty$ is
infinitely divisible, then it satisfies (\ref{indep}).

\section{Size biasing and Compound Poisson}\label{section compound}

Let us now return to our main theme, determining for which distributions
we have (\ref{indep}). We have already seen it is true, trivially, for a scale
multiple of a Poisson random variable. We combine this with
the observation (\ref{noniid}) that to size bias a sum of independent
random variables, just bias a single summand, chosen proportional to its
expectation. Consider a random variable of the form
\begin{equation}\label{poisson sum}
X=\sum_1^n X_j, \ \ \  \mbox{ with }X_j=y_j Z_j, \quad
 Z_j \sim Po(\lambda_j),  \ \
\ \ Z_1,\ldots,Z_n \ \mbox{independent},
\end{equation}
with distinct constants $y_j > 0$.

Since $X$ is a sum of independent variables, we can size bias $X$ by the recipe (\ref{noniid}); pick a summand proportional to its expectation and size bias that one.
We have $EX_j=y_j \lambda_j$ and therefore $a=EX=\sum_j y_j \lambda_j$. Hence, the
probability that we pick summand $j$ to size bias is
$$
\p(I=j)=y_j \lambda_j/a.
$$
But by (\ref{yXpois}), $X_j^*=X_j+y_j$, so that when we pick $X_j$ to bias we
add $y_j$. Hence, to bias $X$ we merely add $y_j$ with probability $y_j \lambda_j/a$,
or, to put it another way $X^*=X+Y$, with $X,Y$ independent and
\begin{equation}\label{discrete-y}
\p(Y=y_j) = y_j \lambda_j / a.
\end{equation}
In summary, $X$ of the form (\ref{poisson sum}) can be size biased by adding an independent, nonnegative increment.
It will turn out that we have now nearly found all solutions of (\ref{indep}),
which will be obtained by taking limits of variables type (\ref{poisson sum}) and adding a nonnegative constant.

\vskip 1pc

Sums of the form (\ref{poisson sum}) are  said to have a
compound Poisson distribution --- of finite type.
Compound Poisson variables in general are
obtained by a method which at first glance looks unrelated to
(\ref{poisson sum}), considering the sum $S_N$
formed by adding a Poisson number $N$ of iid  summands $A_1,A_2,\ldots$ from
{\em any} distribution,
i.e.~taking
\begin{equation}\label{general comp pois}
X \indist S_N := A_1+\cdots + A_N, \ \ \ N  \sim Po(\lambda),
\end{equation}
where $N,A_1,A_2,\ldots$ are independent and $A_1,A_2,\ldots$
identically distributed. The notation $S_N := A_1+\ldots + A_N$ reflects that
of a random walk,  $S_n := A_1 + \cdots A_n$ for $n=0,1,2,\ldots$, with $S_0 =0$.

To fit the sum (\ref{poisson sum}) into the form
(\ref{general comp pois}),
let $\lambda=\sum  \lambda_j$ and  let $A,A_1,A_2,\ldots$  be iid with
\begin{equation}
\label{Afintyp}
\p(A=y_j) = \lambda_j/\lambda.
\end{equation}
The claim is that
$X$ as specified by (\ref{poisson sum}) has the same distribution
as the sum $S_N$.
Checking this claim is most easily done with
generating functions, such as characteristic functions, discussed in the next section.
Nevertheless, it is an enjoyable
exercise for the special case where $a_1,\ldots,a_n$ are mutually
irrational, i.e.~linearly independent over the rationals, so that with $k_1,\ldots,k_n$ integers, the sum $\sum_1^n k_j a_j$ determines the values of the $k_j$.

Note that the distribution (\ref{Afintyp}) of the summands $A_i$
 is {\em different} from the distribution in (\ref{discrete-y}).
In fact, $A^* \indist  Y$, which can be checked using (\ref{Afintyp}) together with (\ref{sizebias-f}), and comparing the result to (\ref{discrete-y}):
\begin{equation}
 \p(A^* = y_j) = {y_j \p(A= y_j) \over \e A }
= {y_j \lambda_j / \lambda \over \sum_k y_k \lambda_k / \lambda }=
{y_j \lambda_j  \over \sum_k y_k \lambda_k  \ }
=\p(Y=y_j).
\end{equation}
Thus the result that for the compound Poisson of finite type,
$X^*=X+Y$ can be expressed as
\begin{equation}\label{one more term}
(S_N)^* = S_N + Y \indist A_1+\cdots + A_N + A^*,
\end{equation}
with all summands independent.
Note how this contrasts with the recipe (\ref{iid}) for size biasing
a sum with a deterministic number of terms $n$;  in the case of (\ref{iid}) the
biased sum and the original sum have the same number of  terms, but in (\ref{one more term}) the
biased sum has one more term than the original sum.

If we want to size bias a general compound Poisson random variable $S_N$,
there must be some restrictions on the distribution for the iid summands $A_i$
in (\ref{general comp pois}).
First, since $\lambda>0$, and (using the independence of $N$ and $S_0,S_1,\ldots$,)
$\e S_N = \e N \ \e A = \lambda \e A$,
the condition that $\e S_N \in (0,\infty)$ is equivalent to
 $\e A \in (0,\infty)$.
The condition that $S_N \geq 0$ is equivalent to the condition $A_i \geq 0$.
We {\em choose} the additional requirement that $A$ be {\em strictly}
 positive, which is convenient since it enables the simple computation $\p(S_N=0)
=\p(N=0)=e^{-\lambda}$.
There is no loss of generality in this added restriction, for if $p := \p(A=0)>0$,
then with $M$ being Poisson with parameter $\e M = \lambda (1-p)$, and $B_1,B_2,\ldots$
iid and independent of $M$,  with $\p(B \in I)=\p(A \in I)/(1-p)$
 for $I \subset (0,\infty)$
(using $p<1$ since $\e A > 0$),
we have $A_1+\cdots+A_N \indist B_1+\cdots + B_M$,
so that $S_N$ can be represented as
a compound Poisson with {\em strictly} positive summands.

\section{Compound Poisson vs. Infinitely Divisible}\label{sect mu nu}
We now have two ways to produce solutions of (\ref{indep}),
infinitely divisible distributions from section \ref{sect inf div}, and finite type
compound Poisson distributions, from (\ref{poisson sum}),
so naturally the next question is:  how are these related?

The finite type compound Poisson random variable in (\ref{poisson sum})
can be extended by adding in a nonnegative constant $c$, to get $X$ of the form
\begin{equation}\label{poisson sum c}
 X=c+\sum_1^n X_j, \ \ \  \mbox{ with }X_j=y_j Z_j, \quad
 Z_j \sim Po(\lambda_j),  \ \
\ \ Z_1,\ldots,Z_n \ \mbox{independent},
\end{equation}
with $c \geq 0$ and distinct constants $y_j > 0$.  In case $c>0$, this random variable is
not called compound Poisson.
Every random variable of the form (\ref{poisson sum c})
is infinitely divisible -- for the $X_i^{(m)}$ [for $i=1$ to $m$ in (\ref{inf div})]
simply take a sum of the same form, but
with $c$ replaced by $c/m$ and $\lambda_j$ replaced by $\lambda_j/m$.

The following two facts are {\em not}  meant to be obvious.
By taking  distributional limits of  sums of the form
(\ref{poisson sum c}), and
requiring that  $ c+ \sum_1^n y_j \lambda_j$ stays bounded as
$n \rightarrow \infty$, one gets all  the non-negative
infinitely divisible distributions with finite mean.
By also requiring  that $c=0$ and $\sum_1^n \lambda_j$
stays bounded,
the limits are all the
finite mean,
non-negative compound Poisson distributions.

To proceed, we calculate the characteristic function for the
distribution in (\ref{poisson sum}).  First, if $X$ is Po($\lambda$), then
$$
\phi_X(u) :=  \e e^{iuX}= \sum_{k \geq 0} e^{iuk} \p(X=k) =
\sum e^{iuk} e^{-\lambda} \frac{\lambda^k}{k!} =
e^{-\lambda} \sum \frac{(\lambda e^{iu})^k}{k!} =\exp(\lambda(e^{iu}-1)).$$
For a scalar multiple $yX$ of a Poisson random variable,
$$
\phi_{yX}(u) = \e e^{iu(yX)} = \e e^{i(yu)X} = \phi_X(yu) = \exp(\lambda(e^{iuy}-1)).
$$
Thus the summand $X_j = y_j Z_j$ in (\ref{poisson sum}) has
characteristic function $\exp(\lambda_j(e^{iuy_j}-1))$,
and hence the sum   $X$ has characteristic function
$$
\phi_X(u) = \prod_{j=1}^n \exp\left(\lambda_j(e^{iuy_j}-1)\right)
= \exp\left( \sum_{j=1}^n \lambda_j(e^{iuy_j}-1)\right).
$$
To prepare for taking limits, we write this as
\begin{equation}\label{finite type}
\phi_X(u) =  \exp\left( \sum_1^n \lambda_j(e^{iuy_j}-1)\right)   =
\exp\left( \ \int_{(0,\infty)} (e^{iuy}-1) \ \mu(dy)\right),
\end{equation}
where $\mu$ is the measure on $(0,\infty)$
which places mass $\lambda_j$ at location
$y_j$.  The total mass of $\mu$ is
$\int 1 \ \mu(dy) = \sum_1^n \lambda_j$, which we will denote by $\lambda$,
 and the
first moment of $\mu$ is $\int y \ \mu(dy) = \sum_1^n y_j \lambda_j$, which
happens to equal $a := \e X$.

Allowing the addition of a constant $c$, the random variable
of the form (\ref{poisson sum c}) has characteristic function $\phi_X$ whose logarithm
has the form
$\log \phi_X(u) = iuc +{ \int_{(0,\infty)}  (e^{iuy}-1) \ \mu(dy) }$, where
$\mu$ is a measure whose support consists of a finite number of points.
The finite mean, not identically zero distributional limits of such random variables
yield all of the finite mean, nonnegative, not identically zero, infinitely divisible distributions.
A random variable $X$ with such a distribution has characteristic function $\phi_X$
with
\begin{equation}
\label{comppois mu}
\log \phi_X(u) = iuc +{ \int_{(0,\infty)}  (e^{iuy}-1) \ \mu(dy) },
\end{equation}
where  $c \geq 0$,
$\mu$ is any nonnegative measure on $(0,\infty)$ such that
$\int y \ \mu(dy) < \infty$, and  not both $c$ and $\mu$ are zero.

Which of the distributions above are compound Poisson? The compound
Poisson variables are the ones in which $c=0$ and $\lambda :=\mu( (0,\infty) ) < \infty$.
With this choice of $\lambda$ for the parameter of $N$ , we have $X \indist S_N$ as
in (\ref{general comp pois}), with the distribution of $A$ given by $\mu/\lambda$, i.e.~$\p(A \in dy) = \mu(dy)/\lambda$.
To check, note that for $X=S_N$, $\phi_X(u) := \e e^{iuS_N}$
$ = \sum_{n\geq 0} \p(N=n) e^{iuS_n}$
$= e^{-\lambda} \sum_{n\geq 0} \lambda^n/n! (\phi_A(u))^n$
$=\exp( \lambda(\phi_A(u) -1)$
$=\exp( \lambda(\int_{(0,\infty)}  (e^{iuy}-1) \ \p(A \in dy)$,
which is exactly (\ref{comppois mu}) with $c=0$.

In our context, it is easy to
tell whether or not a given infinitely divisible random variable is
also compound Poisson --- it is if and only if   $\p(X=0)$ is
strictly positive --- corresponding to $\p(N=0)=e^{-\lambda}$
and $e^{-\infty}=0$.
Among the examples of infinitely divisible random
variables in section \ref{divisible examples}, the only compound Poisson
examples are the Geometric and Negative binomial family, and the distributions
related to Buchstab's function.

In (\ref{discrete-y}) --- specifying the
distribution of the increment $Y$ when $X^*=X+Y$ with $X,Y$ independent ---
the factor $y_j$ on the right hand side
 suggests
another way to write (\ref{comppois mu}).  We multiply and divide by $y$ to get
\begin{equation}\label{comppois}
\log \phi_X(u) = \int_{[0,\infty)} \frac{e^{iuy}-1}{y} \ \nu(dy),
\end{equation}
where $\nu$ is any nonnegative measure on $[0,\infty)$ with total mass
$\nu( \ [0,\infty) \ ) \in (0,\infty)$.
The measure $\nu$ on $[0,\infty)$ is related to $c$ and
$\mu$ by $\nu( \{ 0 \} ) =c$ and for $y>0$, $\ \nu(dy)=y \ \mu(dy)$; we follow the
natural
convention that $(e^{iuy}-1)/y$ for $y=0$  is interpreted as $iu$.
The measure $\nu/a$ is a probability measure on $[0,\infty)$ because
\begin{equation}\label{a from nu}
a :=  \e X = -i(d\log \phi_X(u)/du)\vert_{u=0}=c + \int_{(0,\infty)} y \ \mu(dy) =c+  \nu( (0,\infty) ) = \nu( [0,\infty) ).
\end{equation}
We believe it is proper to refer to either
(\ref{comppois mu}) or (\ref{comppois}) as
a L\'evy representation, and to refer to either $\mu$ or $\nu$ as the \levy measure,
in honor of Paul L\'evy; indeed when restriction that $X$ be nonnegative is dropped, there
are still more forms for the \levy representation,
see e.g.~\cite{feller2} Feller volume II, chapter XVII.

It is not hard to see that any distribution specified by (\ref{comppois})
satisfies (\ref{indep}), by
the following calculation with characteristic functions.
Note that for $g(x):= e^{iux}$, the characterization (\ref{sizebias})
of the distribution of $X^*$ directly gives the characteristic
function $\phi^*$ of $X^*$ as
$\phi^*(u) :=  \e e^{iuX^*} $ $= \e (X e^{iuX})/a$.
For any $X \geq 0$ with finite mean, by an application
of the dominated convergence theorem, if $\phi(u) :=  \e e^{iuX}$
then $\phi '(u) =  \e (iX e^{iuX})$.
Thus for any $X \geq 0$ with $0 < a= \e X < \infty$,
\begin{equation}
\label{phistar}
 \phi^*(u) = \frac{1}{ia} \ \phi'(u).
  \end{equation}
Now if $X$ has characteristic function $\phi$ given by (\ref{comppois}),
 again using dominated convergence,
\begin{equation}\label{check}
\phi'(u)=\phi(u) \  \left(ic +  \int_{(0,\infty)} iy e^{iuy}  \  \mu(dy) \right)
= ia \ \phi(u) \ \int_{[0,\infty)} e^{iuy} \\ \nu(dy)/a.
  \end{equation}
Taking the probability measure $\nu/a$ as the distribution of $Y$, and writing $  \eta$ for
the characteristic function of $Y$, (\ref{check}) says that
$\phi'(u) = ia  \ \phi(u) \   \eta(u)$.  Combined with (\ref{phistar}), we have
\begin{equation}\label{direct}
\phi^* = \phi \   \eta.
\end{equation}
Thus $X^*=X+Y$, with $X$ and $Y$ independent and
${\mathcal L}(Y)=\nu/a$.

For the compound Poisson case in general,
in which the distribution of $A$ is
$\mu/\lambda$,
we have $Y\indist A^*$ because $a := \e X = \lambda \ \e A$ and
\begin{equation}\label{A gen type}
\p(A^* \in dy) = y \p(A \in dy)/\e A ={ y (\mu( dy)/\lambda) \over a/\lambda} =
\nu(dy)/a = \p(Y \in dy),
\end{equation}
which can be compared discrete version (\ref{Afintyp}).  Thus the computation (\ref{direct}) shows, for the compound Poisson case, that
(\ref{one more term}) holds.

\section{Main Result}

Our main result is essentially the converse of the computation (\ref{direct}).

\begin{theorem}\label{thm}\label{Main Result}(Steutel 1973 \cite{steutel} )
For a  random variable $X \geq 0$ with $a :=  \e X \in (0,\infty)$,
the following three statements are equivalent.

i) There exists a coupling with $X^*=X+Y, Y \geq 0$, and $X,Y$ independent.

ii) The distribution of $X$ is infinitely divisible,

iii) The characteristic function of $X$ has the L\'evy representation (\ref{comppois}).

Furthermore, when any of these statements hold,
the \levy measure $\nu$ in (\ref{comppois}) equals $a$
times the distribution of $Y$.
\end{theorem}

\proof
We  have proved that ii) implies i), in section \ref{sect inf div}.
We have proved that iii) implies i), in the argument ending with (\ref{direct}), which also
shows that given $\nu$ in the \levy representation (\ref{comppois}), the increment
$Y$ in i) has distribution $\nu/a$.

The equivalence of ii) and iii) is a standard textbook topic --- with the argument for
iii) implies ii) being simply that $X$ with a given $\nu$ is the sum of $n$ iid terms each having the \levy representation (\ref{comppois}) with $\nu/n$ playing the role of $\nu$.

Now to prove that i) implies ii), we assume that i) holds.
The characteristic function $\phi^*$ of $X^*$ has the
form $\phi^* = \phi \,   \eta$, where $\phi$ and $  \eta$ are
the characteristic functions of $X$ and $Y$, so that
$  \eta(u) =  \e e^{iuY} = \int_{[0,\infty)} e^{iuy} \ \p(Y \in dy)$.
Combining this with (\ref{phistar}) we have
$$
\frac{1}{ia} \ \phi'(u) = \phi^*(u) = \phi(u) \   \eta(u)
$$
so that $(\log \phi(u))'=i a \   \eta(u)$.  Since $\log \phi(0)=0$,
 \beas
\log \phi(u) &=& {i a \int _{s\in [0,u)}   \eta(s) \ ds } \\
        &=& {i a \int _{s\in [0,u)} \int_{y \in [0,\infty)} e^{isy}\ \p(Y \in dy)\ ds} \\
        &=& {i a \int _ {y \in [0,\infty)} \int_{s\in [0,u)} e^{isy} \ ds\ \p(Y \in dy)} \\
        &=& {i a \left(
\int _ {y \in (0,\infty)} \int_{s\in [0,u)}  e^{isy} \ ds\ \p(Y \in dy) \
+ u \p(Y=0) \right) }\\
        &=& { a \int _{y \in [0,\infty)}  \frac{e^{iuy} -1}{y} \  \p(Y \in dy) }.
  \enas
This is the same as the representation (\ref{comppois}), with
$\nu = a {\mathcal L}(Y)$ for the random variable $Y$ given in i).

\qed

Observe that $\nu$ is an arbitrary probability distribution on $[0,\infty)$,
i.e. $\nu \in $Pr($[0,\infty))$,
and the choice of $a \in (0,\infty)$ is also arbitrary.  Thus there is
a one-to-one correspondence between
the Cartesian product $ $Pr($[0,\infty)) \times (0,\infty)$ and
the set of the nonnegative, infinitely divisible distributions with
finite, strictly positive mean.

\section{A consequence of  $X^*=X+Y$ with independence}

To paraphrase the result of Theorem \ref{Main Result},
for a nonnegative
random variable $X$ with $a := \e X \in (0,\infty)$, it is possible
to find a coupling with  $X^*=X+Y$, $Y\geq 0$
and $X,Y$ independent {\em if and only if} the distribution of $X$ is the infinitely
divisible distribution with L\'evy representation (\ref{comppois})
governed by the finite measure $\nu$ equal to $a$ times the distribution of
$Y$.
Thus we know an explicit, albeit complicated, relation between the distributions of
$X$ and $Y$.  It is worth seeing
how (\ref{indep}) directly gives a simple relation between the densities of $X$ and $Y$,
if these densities exist.

In the discrete case, if $X$ has a mass function $f_X$
and if (\ref{comppois}) holds,
then $Y$ must mass function, $f_Y$, and by (\ref{indep}),
$f_{X^*}$ is the convolution
of $f_X$ and $f_Y$: \ $f_{X^*}(x) = \sum_y f_X(x-y) f_Y(y)$.  Combined with
(\ref{sizebias-f}), this says that for all $x>0$,
$$
f_X(x) = \frac{a}{x} \sum_y f_X(x-y) f_Y(y).
$$
Likewise, in the continuous case, if  $X$ has density $f_X$
(i.e. if for all bounded $g$,
$\e g(X) = \int g(x)f_X(x) \ dx$,)
and if (\ref{comppois}) holds,   {\em and if further}
$Y$ has a density $f_Y$, then by (\ref{indep}),
$f_{X^*}$ is the convolution
of $f_X$ and $f_Y$: \ $f_{X^*}(x) = \int_y f_X(x-y) f_Y(y)$.  Combined with
(\ref{sizebias-f}), this says that for all $x>0$,

\begin{equation}\label{cont convolve}
f_X(x) = \frac{a}{x} \int_y f_X(x-y) f_Y(y) \ dy .
\end{equation}

\section{Historical remark}
For all intents and purposes, Theorem \ref{thm} is due to
Steutel \cite{steutel}.  The way he states his result is sufficiently
different from our Theorem \ref{thm} that for comparison,
we quote verbatim from \cite{steutel}, p. 136:

\noindent {\bf Theorem 5.3.}  A d.f. $F$ on $[0,\infty)$ is
infinitely divisible iff it satisfies
$$
(5.6) \ \ \ \
\int_0^x u \ dF(u) \ = \ \int_0^x F(x-u) \ dK(u),
$$
where $K$ is non-decreasing.

Observe that Steutel's result is actually more general than Theorem \ref{thm},
since that latter only deals with nonnegative infinitely divisible random
variables with {\em finite mean}.
The explicit connection between the independent increment for size biasing,
and the \levy representation, is made in \cite{harn and steutel}, along with
further connections between renewal theory and independent increments.

\section{The product rule for size biasing}

We have seen that for independent, nonnegative random variables
$X_1, \ldots, X_n$, the sum $X=X_1+X_2\cdots+X_n$ can be size biased
by picking a single summand at random with probability proportional to its expectation,
and replacing it with one from its size biased distribution.
Is there a comparable procedure for the product $W= X_1 X_2 \cdots X_n$?
Would it involve size-biasing a single factor?

Let $a_i = \e X_i \in (0,\infty)$, let $F_i$ be the distribution function of
$X_i$, and let
$F_i^*$ be the distribution function of $X_i^*$, so that
$dF_i^*(x)=x\ dF_i(x)/a_i$.  Let $X_1^*,\ldots,X_n^*$ be independent.
By (\ref{sizebias}) with $a := \e W = a_1 a_2 \cdots a_n$,
for all bounded functions $g$,
\beas
Eg(W^*) &=& E\left(W g(W)\right)/ (a_1 a_2 \cdots a_n)\\
        &=& \int \cdots \int x_1\cdots x_n \ \ g(x_1 x_2 \cdots x_n) \ dF_1(x_1)\cdots dF_n(x_n)/(a_1\cdots a_n)\\
        &=& \int \cdots \int g(x_1 x_2 \cdots x_n) \ (x_1\ dF_1(x_1)/a_1)
\cdots (x_n \ dF_n(x_n)/a_n)\\
        &=& \int \cdots \int g(x_1 x_2 \cdots x_n) \  dF_1^*(x_1)\cdots \ dF_n^*(x_n) \\
        &=&  Eg(X_1^* \cdots X_n^*),
\enas
and so
$$
W^* \indist X_1^*\cdots X_n^*.
$$
We have shown that to size bias a product of independent variables, one must  size bias
every factor making up the product, very much unlike what happens for a sum,
where only one term is size biased!

\section{Size biasing the lognormal distribution}\label{sect lognormal}
The lognormal distribution is often used in financial mathematics to model
prices, salaries, or values.
A variable $L$ with the lognormal distribution
is obtained by exponentiating a normal variable.
We follow the convention that $Z$ denotes a
standard normal, with $\e Z=0,$
var $Z =1$, so that $L=e^Z$ represents a standard lognormal.  With
constants $\sigma>0, \mu \in \BR$,
$\sigma Z + \mu$ represents the general normal, and $L=e^{\sigma Z + \mu}$
represents the general lognormal.
As the lognormal is non-negative and has finite mean, it
can be size biased to form $L^*$.

One way to guess the identity of $L^*$ is to use the method of moments.
For the standard  case $L=e^Z$, for any
real $t$, calculation gives $\e e^{tZ} = \exp(t^2/2)$.  Taking
$t=1$ shows that $\e L= \sqrt{e}$, and more generally, for
$n=1,2,\ldots$, $\e L^n = \e e^{nZ} = \exp(n^2/2).$
Using relation (\ref{moment shift}), the moment-shift for size biasing,
we have $\e (L^*)^n = \e L^{n+1}/\e L$ $=\exp((n+1)^2/2 - 1/2)$
$=\exp(n^2/2 + n)$
$=e^n \e L^n = \e (eL)^n$.  Clearly we should guess
that $L^* \indist eL$, but we must be cautions, as
the most famous example of a distribution which has moments
of all orders but which is not determined
by them is the lognormal; for other such distributions
related to the normal, see \cite{slud}.

We present a rigorous
method for finding the distribution of $L^*$,
based on the size biasing product rule of the previous section;
as an exercise the reader might try to verify our conclusion
(\ref{sb lognormal})
by working out the densities for lognormal distributions, and
using the relation (\ref{sizebias-f}).

We begin with the case $\mu=0, \sigma>0$.
Let $C_i$ be independent variables taking the values $1$ or $-1$ with equal
probability. These variables have mean zero and variance one, and by the
central limit theorem, we know that
$$
\frac{1}{\sqrt{n}} \sum_{i=1}^n \sigma C_i \todist \sigma Z.
$$
Hence, we must have
$$
W=\prod_{i=1}^n \exp(\frac{1}{\sqrt{n}}\sigma C_i)=
\exp(\frac{1}{\sqrt{n}} \sum_{i=1}^n \sigma C_i) \todist \exp(\sigma Z)=L,
$$
a lognormal, and thus $W^* \todist L^*$.
Write  $X_i:=\exp(\sigma C_i/\sqrt{n})$,
so that $W=X_1\cdots X_n$ with independent factors, and
by the product rule,
$
W^*=  X_1^*\cdots X_n^*.
$
The variables $X_i$ take on the values
$q=e^{-\sigma/ \sqrt{n}}$ and $p=e^{\sigma/ \sqrt{n}}$ with equal probability,
and so $X_i^*$ take on these same values, but with probabilities
$q/(p+q)$ and $p/(p+q)$ respectively. Let's say that $B_n$
of the $X_i^*$ take the value $p$, so that $n-B_n$ of the $X_i^*$
take the value $q$.
Using $B_n$, we can write
$$
W^*=p^{B_n}q^{n-B_n}=e^{\sigma (2B_n-n)/\sqrt{n}}.
$$
Since $B_n$  counts the number of ``successes'' in $n$ independent
trials, with success probability $p/(p+q)$, $B_n$ is distributed binomial$(n,p/(p+q))$.
As $n \rightarrow \infty$, the central limit theorem gives that $B_n$ has an
approximate normal distribution.
Doing a second order Taylor expansion of $e^x$ around zero, and applying
it at $x = \pm \sigma / \sqrt{n}$, we find that
$p/(p+q) = 1/2 + \sigma /(2 \sqrt{n})+O(1/n)$, so that
$B_n$ is approximately normal, with mean
$np/(p+q)=(1/2)(n+ \sigma \sqrt{n})+O(1)$ and variance
$npq/(p+q)^2=n/4+O(1/n^{3/2})$.  Hence
$$
\frac{1}{\sqrt n}(2B_n-n) \todist   Z+\sigma  \quad \mbox{as $n \rightarrow \infty$}
$$
and therefore
$$
W^* \todist  e^{\sigma (Z+\sigma)}.
$$
Since $W^* \todist L^* = (e^{\sigma Z})^*$, we have shown that $(e^{\sigma Z})^*
\indist e^{\sigma (Z+\sigma)}$.
For the case where $L=e^{\sigma Z + \mu}$,
the scaling relation (\ref{scaling}) yields
the formula for size biasing the lognormal in general:
\begin{equation}\label{sb lognormal}
(e^{\sigma Z + \mu})^*=e^{\sigma (Z+\sigma) + \mu}.
\end{equation}

\section{Examples}
In light of Theorem \ref{thm}, for a nonnegative random variable $X$
 with finite, strictly positive mean,
being able to satisfy
 $X^*=X+Y$ with independence and $Y \geq 0$  is
equivalent
to being infinitely divisible.
We give examples of size biasing, first with examples that are not infinitely
divisible, then with examples that are.

\subsection{ Examples of size biasing without an independent increment}
\label{exnoindep}

Both  examples 1 and 2 below involve
bounded, nonnegative random variables.
Observe that in general, the distributions of  $X$ and $X^*$ have the same
support, except that always $\p(X^*=0)=0$.  This immediately implies that if $X$ is
bounded but not constant, then it cannot satisfy (\ref{indep}).

\vskip 1pc \noindent {\bf Example 1.  Bernoulli and binomial}
\label{bab}

Let $B_i$ be Bernoulli with parameter $p \in (0,1]$, i.e.~$B_i$ takes the value 1 with probability $p$,
and the value 0 with
probability $1-p$.  Clearly $B_i^*=1$, since $\p(B_1^*=1)=1\p(B_1=1)/ \e B_1=1$.
If $B_1,B_2,\ldots$ are independent, and
$S_n=B_1+\cdots+B_n$
we say that $S_n \sim$ binomial $(n,p)$.
We size bias $S_n$ by size biasing a single summand, so
$S_n^* \indist S_{n-1}+1$,
which cannot be expressed as $S_n+Y$ with $S_n,Y$ independent!

Note that letting
$n \rightarrow \infty$ and $np \rightarrow \lambda$ in the relation
$S_n^* \indist S_{n-1}+1$
gives another proof that $X^* \indist X+1$ when $X \sim Po(\lambda)$,
because both $S_{n-1} \todist X$ and $S_n \todist X$.
Here we have a family of examples without
independence, whose limit is the basic example with independence.

\vskip 1pc \noindent {\bf Example 2.  Uniform and Beta}
\label{uab}

The Beta distribution on (0,1), with parameters $a,b>0$
is specified by saying that its has a density on (0,1),
proportional to $(1-x)^{a-1}x^{b-1}$.
The uniform distribution on (0,1) is the special case $a=b=1$
of this Beta family.  Using (\ref{sizebias-f}), if
$X \sim$ Beta$(a,b)$, then $X^* \sim $ Beta$(a,b+1)$.

\vskip 1pc
There are many families of distributions for which size biasing simply changes the
parameters; our examples are the Beta family in example 2,
the negative binomial family in example 4, the Gamma family in example 5,  and
the lognormal family in example 6.
In these families, either all members satisfy (\ref{indep}), or else
none do.  Thus it might be tempting to guess that infinite divisibility is a property
preserved by size biasing, but it ain't so.

\noindent {\bf Example 3.   $X=1+W$ where $W$ is Poisson}
\label{1plusP}

We have $X^*$ is a mixture of $X+0$ and $X+1$, using (\ref{noniid}) with
$X_1=1$, $X_1^*=X_1+0$ and $X_2=W$, $X_2^*=W+1$.  That is,
$X^*$ is a mixture of $1+W$, with weight $1/(1+\lambda)$, and $2+W$, with
weight $\lambda/(1+\lambda)$.  Elementary calculation shows that it is not
possible to have $X^*=X+Y$ with $X,Y$ independent and $Y\geq 0$.
Therefore $X$ is not infinitely divisible.

Since $X=W^*$, we have an example in which $W$ is infinitely divisible, but
$W^*$ is not.

\vskip 1pc
\subsection{\bf Examples of $X^*=X+Y$ with independence}\label{divisible examples}

By Theorem \ref{thm}, when $X$ satisfies $X^* \indist X+Y$ with $X,Y$ independent
and $Y \geq 0$, the distribution of $X$ is determined by the distribution of $Y$ together
with a choice for the constant $a \in (0,\infty)$ to serve as $\e X$.
Thus all our examples below, organized by a choice of $Y$, come
in one parameter families indexed by $a$ --- or if more convenient, by
something proportional to $a$;
in these families, $X$ varies and $Y$ stays constant!

\vskip 1pc \noindent {\bf Example 4. $Y$ is 1+geometric.  $X$ is  geometric or negative binomial}
\label{1plusG}

4a) The natural starting point is that you are given the geometric distribution:
$\p(X = j) = (1-q)q^j$  for
$j \geq 0$, with $0<q<1$, and you want to discover whether or not it is
infinitely divisible. Calculating the characteristic function,
$\phi(u) = \sum_{k \geq 0} e^{iuk}(1-q)q^k = (1-q)/(1-qe^{iu})$,
so $\log \phi(u) = \log(1-q) - \log(1-qe^{iu})$
$=-\sum_{j \geq 1} q^j/j + \sum_{j \geq 1}(q^je^{iuj}/j)$
$=\sum_{j \geq 1} ((e^{iuj}-1)/j) \ q^j$.

Thus the geometric distribution has a \levy representation in which
$\nu$ has mass $q^j$  at
$j=1,2,\ldots$, so we have verified that the geometric distribution is infinitely divisible.
The total mass $a$ of $\nu$ is $a= q+q^2+\cdots = q/(1-q)$; and
this agrees with the recipe $a = \e X$.  Since $\p(Y=j)
= \nu(\{j\})/a = (1-q)q^{j-1}$ for
$j=1,2,\ldots$, we have $Y \indist 1+X$.
Thus $X^*=X+ Y$ with $X,Y$
independent and $Y \indist X+1$.

4b)
Multiplying the \levy measure $\nu$ by $t>0$ yields the
general case of the negative binomial distribution, $X \sim $ negative binomial$(t,q)$.
The case $t=1$ is
the geometric distribution.
We still have $X^*=X+Y$ with $X,Y$
independent, and $Y \sim $ geometric$(q)+1$. Note that for integer $t$ we can
verify our calculation in another way, as in this case $X$ is
the sum of $t$ independent geometric($q$) variables $X_i$. By (\ref{iid}), we can
size bias $X$ by size biasing a single geometric term, which
is the same as adding an independent $Y$ with distribution, again, 1 + geometric($q$).

\vskip 1pc \noindent {\bf Example 5.  $Y$ is exponential.  $X$ is exponential or Gamma}
\label{YisE}

5a)
Let  $X$ be exponentially distributed with $ \e X=1/\alpha$, i.e. $\p(X>t)=e^{-\alpha t}$ for
$t > 0$.  As we saw in section \ref{sec:renewal} for the case $\alpha=1$,
 $X^*=X+ Y$ with $X,Y$
independent and $Y \indist X$. The case with general $\alpha >0$  is simply
a spatial scaling of the mean one, ``standard'' case.  The L\'evy measures $\nu$ is
simply $\e X = 1/\alpha$ times the common distribution of $X$ and $Y$, with
$\nu(dy) = e^{-\alpha y}\  dy$

5b) Multiplying the \levy measure $\nu$ by $t>0$ yields the
general case of the Gamma distribution, $X \sim $ Gamma$(\alpha,t)$.
The name comes from that fact that $X$ has density
$f(x)= (\alpha^t/\Gamma(t)) \ x^{t-1} e^{-\alpha x}$ on
$(0,\infty)$.
The special case $t=1$ is the exponential distribution, and more generally the
case $t=n$ can be realized as $X = X_1+\cdots+X_n$ where the $X_i$ are iid,
exponentially distributed with $\e X_i = 1/\alpha$.
We have  $X^*=X+ Y$ with $X,Y$
independent and $Y$ is exponentially distributed with
 mean ($1/\alpha)$, so that
$X^* \sim $ Gamma$(\alpha,t+1)$. The  \levy measure here is $\nu$
with $\nu(dy) =t  e^{-\alpha y}\  dy$; so the corresponding $\mu$ has
$\mu(dy) = t e^{-\alpha y} / y \ dy$.  This form of $\mu$ is
known as the Moran or Gamma subordinator; see e.g.
\cite{kingman}. As in example 4b), for integer $t$ we can
verify our calculation by noting that $X$ is
the sum of $t$ independent exponential, mean ($1/\alpha$) variables, and that by
(\ref{iid}), when size biasing we will get the same $Y$ added on to the sum as the $Y$ which
appears when size biasing any summand.

\vskip 1pc \noindent {\bf Example 6.  $Y$ is ??,  \ $X$ is lognormal}

As mentioned in Section (\ref{sect lognormal}), we say that $X$
is lognormal when $X = e^{\sigma Z + \mu}$ where $Z$ is a standard normal variable.
The proof that the lognormal is infinitely divisible, first given by
Thorin \cite{thorin}, remains difficult;
there is an excellent book by Bondesson \cite{bondesson}
for further study.  Consider even the standard case, $X=e^Z$, so that by equation (\ref{sect lognormal}), $X^* \indist e^{Z+1} \indist e X$.  The result of Thorin that this $X$ is
infinitely divisible is thus equivalent, by Theorem \ref{thm}, to the claim that there
exists a distribution for $Y \geq 0$ such that with $X$ and $Y$ independent, $X+Y \indist eX$.
Also, by Theorem \ref{thm} with $a=\e X = \sqrt{e}$,
 the distribution of $Y$ is exactly $1/\sqrt{e}$ times the
\levy measure $\nu$.  However, there does
not seem to exist yet a simplified expression for this distribution!

Since the lognormal $X=e^Z$ satisfies $X^*=eX=X+(1-e)X$, it provides a simple illustration
of our remarks in the paragraph following (\ref{xplusy}),
that the relation $X^*=X+Y, Y \geq 0$ does not determine the distribution of $Y$ without
the further stipulation that $X$ and $Y$ be independent.  Note also that in
 $X^*=X+(1-e)X$, the increment $(1-e)X$ is a monotone function of $X$, so this is an example of the coupling
using the quantile transformation.

\vskip 1pc \noindent {\bf Example 7.  $Y$ is uniform on an interval $(\beta,\gamma)$, with $0 \leq \beta < \gamma < \infty$.}
By scaling space (dividing by $\gamma$) we can assume without loss of
generality that $\gamma=1$.  This still allows two qualitatively distinct cases,
depending on whether $\beta=0$ or $\beta>0$.

\vskip 1pc \noindent {\bf Example 7a.  $\beta=0$: \
Dickman's function and its convolution powers}.

With $a=\e X \in (0,\infty)$,
this example is specified by (\ref{comppois}) with
$\nu$ being
$a$ times the uniform distribution  on (0,1), so that $\mu(dx)=a/x \ dx$ on $(0,1)$.
The reader must take on faith that $\nu$ having a density,
together with $\mu( \ (0,\infty) \,)= \infty$ so that $\p(X=0)=0$,
implies that
the distribution of $X$ has a density, call it
$g_a$.   Size biasing then gives an interesting differential-difference equation
for this density:  using (\ref{cont convolve}),  for $x>0$,
\begin{equation}\label{dickman int}
  g_a(x) = \frac{a}{x} \int_{y=0}^1 g_a(x-y) \ dy
= \frac{a}{x} \int_{x-1}^x g_a(z) \ dz.
\end{equation}
Multiplying out gives $x g_a(x)=a\int_{x-1}^x g_a(z) \ dz$, and
taking the derivative with respect to $x$ yields
$x g_a'(x)+ g_a(x) = a g_a(x)-a g_a(x-1)$, so that
$g_a'(x) =( \ (a-1) g_a(x) - a g_a(x-1) \ ) / x$, for $x>0$.

For the case
$a=1$ this simplifies to
$g_1'(x)=-g_1(x-1)/x$, which is the same differential-difference
equation that is used to specify Dickman's function $\rho$,
of central importance in number theory; see \cite{tenenbaum}.
The function $\rho$ is characterized by $\rho(x)=1$ for $0 \leq x \leq 1$ and $\rho'(x)=-\rho(x-1)/x$ for $x>0$, and $\rho(x)=0$ for $x<0$, with $\rho$ continuous on $[0,\infty)$, and from the calculation that $\int_0^\infty \rho(x) \ dx = e^\gamma$,
where $\gamma$ is Euler's
constant,  it follows that $g_1(x) = e^{-\gamma} \rho(x)$.
Dickman's function
governs the distribution of the
largest prime factor of a random integer in the following sense:  for fixed $u>0$,
the proportion of integers from 1 to $n$ whose largest prime factor is smaller than
$n^{1/u}$ tends to $\rho(u)$ as $n \rightarrow \infty$.
For example, $\rho(2)$ can be calculated from the differential equation
simply by $\rho(2)= \rho(1) + \int_1^2 \rho'(x) \ dx $
$= 1 + \int_1^2 -\rho(x-1)/x \ dx$
$= 1 + \int_1^2 -1/x \ dx$
$=1-\log 2 \doteq 1-.69314 = .30686$, and
the claim is that $\rho(2)$ gives, for large $n$, the approximate
proportion of integers from 1 to $n$ all of whose prime factors are
at most $\sqrt{n}$.

For general $t>0$ the density $g_t$ is a ``convolution power of Dickman's function,''
see \cite{hensley}.
The size bias treatment of this first appeared in the 1996 version of \cite{abt}, and was subsequently written up in \cite{A98}.

\vskip 1pc \noindent {\bf Example 7b. $\beta > 0$: \  Buchstab's function, integers free of small prime factors.}

For these examples $Y$ is uniform on $(\beta,1)$ for $\beta \in (0,1)$,
with density $1/(1-\beta)$ on $(\beta,1)$.  Therefore
$\nu$ is a multiple of  uniform distribution on $(\beta,1)$, with
density $t$ on $(\beta,1)$ for some constant $t>0$
--- we have $a:=\e X = t(1-\beta)$ --- but
$t$ rather than $a$ is the convenient parameter.
From $\nu(dx) = t \ dx$ on $(\beta,1)$ we get
$\mu(dx) = t/x \ dx$ on $(\beta,1)$, so that the
total mass of $\mu$ is $\lambda = \int_{(\beta,1)} t/x \ dx $
$= t \log(1/\beta)$.  Since $\lambda < \infty$,   $X$ is compound Poisson with
$\p(X=0) = e^{-\lambda} = \beta^t$.

For the case $t=1$, the  distribution of the random variable $X$ is related to
another important function in number theory, Buchstab's function $\omega$;
again see \cite{tenenbaum}.
The relation involves a ``defective density''  --- here $t=1$ and
 $\p(X=0)=\beta >0$ so $X$ does not have a proper density.
Size biasing yields a relation similar to (\ref{cont convolve}),
which leads to a differential-difference equation, which in turn
establishes the relation between the defective density and
Buchstab's function; see \cite{AS}.  The net result is that
for $\beta < a <b < 1$,
$\p( a<X<b)=\int_a^b \omega(x/\beta) \ dx$.
 Buchstab's function $\omega$ is characterized by the properties that it
is continuous on $(1,\infty)$,  $\omega(u)=1/u$ for $u\in [1,2]$, and
$(u \omega(u))'=\omega(u-1)$ for $u>2$.
It governs the distribution of the
smallest  prime factor of a random integer in the sense that for $u>1$,
the proportion of integers form 1 to $n$ whose smallest prime factor is at least $n^{1/u}$
is asymptotic to $ u \omega(u) / \log n$.

\end{document}